\documentclass[12pt]{article}
\usepackage{amsxtra,amssymb,amsthm,amsmath,latexsym}

\theoremstyle{plain}

\def\R{{\mathbb R}}

\def\oH{{\overset{\circ}{H}}}
\def\oH1{{\overset{\circ}{H}\kern-.02in{}^1}}

\def\bee{\begin{equation*}}
\def\eee{\end{equation*}}
\def\be{\begin{equation}}
\def\ee{\end{equation}}

\begin{document}

\title{On a hyper-singular equation
}

\author{Alexander G. Ramm\\
 Department  of Mathematics, Kansas State University, \\
 Manhattan, KS 66506, USA\\
ramm@math.ksu.edu\\
http://www.math.ksu.edu/\,$\sim$\,ramm}

\date{}
\maketitle\thispagestyle{empty}

\begin{abstract}

The equation $v=v_0+\int_0^t(t-s)^{\lambda -1}v(s)ds$ is considered, $\lambda\neq 0,-1,-2...$ and $v_0$ is a smooth function rapidly decaying with all its derivatives.
It is proved that the solution to this equation does exist, is unique and is  smoother than the singular function $t^{-\frac 5 4}$. 

\end{abstract}

\section{Introduction and formulation of the result}\label{S:1}

Let
 \be\label{e1}
v(t)=v_0(t)+\int_0^t (t-s)^{\lambda -1}v(s)ds,
\quad \lambda\neq 0,-1,-2,....
\ee 
where  $v_0$ is a smooth functions rapidly decaying with all its derivatives as $t\to \infty$, $v_0(t)=0$ if $t<0$. The integral in \eqref{e1} diverges in the classical sense.

Our result can be formulated as follows.

{\bf Theorem 1.} {\em The solution to equation \eqref{e1} for $\lambda=-\frac1 4$ does exist, is unique and less singular than 
$t^{-\frac 5 4}$ as $t\to 0$.} 

{\em Proof of Theorem 1.} We define the integral in \eqref{e1} as a convolution of the distribution $t^{\lambda-1}$ and $v$.
The space of the test functions for this distribution is the space $\mathcal{K}:=C^\infty_0(\R_+)$ of compactly supported on $\R_+:=[0, \infty)$ infinitely differentiable functions $\phi(t)$ defined on $\R_+:=[0,\infty)$. The topology on this space is defined by 
countably many norms $\sup_{t\ge 0}t^m|D^p\phi(t)|$.  A sequence $\phi_n(t)$ converges
to $\phi(t)$ in $\mathcal{K}$ if and only if
all the functions $\phi_n(t)$ have compact support on 
an interval $[a,b]$, $a>0$, $b<\infty$ and
$\phi_n$ converges on this interval to $\phi$
in every of the above norms. 
 Let us check that $t^{\lambda-1}:=t_+^{\lambda-1}$
is a distribution on  $\mathcal{K}$ for $\lambda<0$, i.e., a linear bounded functional on $\mathcal{K}$. Let $\phi_n\in \mathcal{K}$ and
$\phi_n \to \phi$ in $\mathcal{K}$. If $\lambda<0$ then
$\max_{t\in [a,b]}t^{\lambda -1}\le a^{\lambda-1}+b^{\lambda-1}$. Thus, 
$$|\int_0^\infty t^{\lambda-1}\phi_n(t)dt|\le
[a^{\lambda-1}+b^{\lambda-1}] \int_0^\infty |\phi_n(t)|dt,$$ where $a>0$  and $b<\infty$. 
Since $\phi_n \to \phi$ in $\mathcal{K}$, we have 
$$\int_0^\infty |\phi_n(t)|dt\to \int_a^b |
\phi|dt.$$ 
So, the integral 
 $\int_0^\infty t^{\lambda-1}\phi(t)dt$ is a bounded linear functional on $\mathcal{K}$ and $t^{\lambda-1}$ is a distribution on the set of the test functions 
 $\mathcal{K}$ for $\lambda\neq 0,-1,-2,....$.

The integral in \eqref{e1} is the convolution $t^{\lambda-1}\star v$. This convolution is 
defined for any distributions on the dual to 
$\mathcal{K}$ space $\mathcal{K}'$. This is done in  \cite{R691}, p.57.  
For another space of the test functions $K=C^\infty_0(\R)$ this is done in 
\cite{GS}, p.135. 

It is  known, see e.g. \cite{BP}, p.39,
that 
 \be\label{e2a} 
L(f\star h)=L(f)L(h),
\ee
 where $L$ is the Laplace transform, and $f,h$ are distributions on $\mathcal{K}$.
 
  Let us calculate $L(t^{\lambda-1})$
using the new variable $s=pt$: 
 \be\label{e2} 
L(t^{\lambda-1})=\int_0^\infty t^{\lambda-1}e^{-pt}dt=\int_0^\infty s^{\lambda-1}e^{-s}ds
p^{-\lambda}=\Gamma(\lambda)p^{-\lambda},\quad
\lambda\neq 0,-1,-2....
\ee
This formula is valid classically for $Re \lambda>0$. By analytic continuation
with respect to $\lambda$ 
it is valid for all complex $\lambda\neq 0,-1,-2,....$.

Applying the Laplace transform to \eqref{e1}
and using formulas \eqref{e2a} and \eqref{e2}, one gets
 \be\label{e3} 
L(v)=L(v_0)+\Gamma(\lambda)p^{-\lambda}L(v).
\ee
Let us assume that $\lambda=-\frac 1 4$ so that $\lambda-1=-\frac 5 4$.
This value appears in the solution to the Navier-Stokes problem in $\R^3$, see \cite{R691}, p.53. If $\lambda=-\frac 1 4$,
then equation \eqref{e3} yields
 \be\label{e4}
 L(v)=\frac{L(v_0)}{1+4\Gamma(3/4)p^{1/4}},
 \ee  
where we have used the relation $\Gamma(-\frac 1 4)=-4\Gamma(3/4)$, which follows from the known formula $\Gamma(z+1)=z\Gamma(z)$ with $z=-\frac 1 4$.

Thus,
 \be\label{e5}
 v=L^{-1}\frac{L(v_0)}{1+4\Gamma(3/4)p^{1/4}}.
 \ee
 So, {\em the solution $v$ does exist and is unique.}   
 
  Moreover, $v$ is not a distribution if $v_0$ is smooth and rapidly decaying when
 $t\to \infty$. This follows from the known results concerning the relation of asymptotic of $L(f)(p)$ and $f(t)$ for $p\to \infty$ and
 $t\to 0$ and for $p\to 0$ and $t\to \infty$,
 see \cite{BP}, p.41.  

Namely, if  $f(t)\sim At^{\nu}$ as $t\to 0$,
then $L(f)(p)\sim A\Gamma(\nu +1)p^{-\nu -1}$
as $p\to \infty$, $\nu\neq -1.-2,....$.
If $f(t)\sim  At^{\nu}$ as $t\to \infty$
then $L(f)(p)\sim  A\Gamma(\nu +1)p^{-\nu -1}$
as $p\to 0$.

 Since $p^{1/4}\to 0$ as $p\to 0$,
the asymptotic of $v(t)$ as $t\to \infty$
is of the same order as that of $v_0$.

As $p\to \infty$ the singularity of $v(t)$ as $t\to 0$ is of the order less than 
that of $t^{-\frac 5 4}$.
For example, assume that $v_0$ is continuous as $t\to 0$.
Then we can take $\nu\ge 0$. Consider the 
worst case $\nu=0$. In this case  $L(v)(p)$ is of the order 
$p^{-1-\frac 1 4}=p^{-\frac 5 4}$. Therefore
$v\sim t^{-\frac 1 4}$ as $t\to 0$. This is 
an integrable singularity.
 Thus, $v$ is less singular as $t\to 0$
 than the distribution $t^{-\frac 5 4}$.

Theorem 1 is proved. \hfill$\Box$

In \cite{R677} another  result, similar to the one in this paper, is proved. 

 In   Zbl 07026037 in a review of paper \cite{R684}
there is an erroneous claim that the proof
in \cite{R684} is incorrect. The reviewer
 erroneously claims that the integral \eqref{e1} diverges and therefore
it is equal to infinity. While this is true classically it is not true in the sense of distributions. Therefore, the claim of the reviewer that the proof in \cite{R684} is not correct is false.
The reviewer
claims that $\Phi_{-\frac 1 4}$ is {\em not} equal to
$\frac {t_+^{-\frac 5 4}}{\Gamma(-\frac 1 4)}$.
This is not true if the space of the test
functions is $\mathcal{K}$ (although it is true if the space of the test functions 
is $K$).

{\bf Concluding remark.} Historically it is 
well known that equation \eqref{e1} can be solved 
explicitly by the Laplace transform if
$\lambda>0$ and the function $1-L(t^{\lambda -1})\neq 0$. To our knowledge, for $\lambda<0$ there were no results concerning the solvability of equation \eqref{e1}. The author 
got interested in \eqref{e1} in the case $\lambda=-\frac 1 4$ in connection with the millennium problem about unique global solvability of the Navier-Stokes problem (NSP)
in $\R^3$ which was solved in \cite{R684}, 
see also \cite{R691} Chapter 5. 
results



\begin{thebibliography}{1000} 

\bibitem{BP}  Yu. Brychkov, A. Prudnikov,
{\bf Integral tranforms of generalized functions}, Nauka, Moskow, 1977  (in Russian)

\bibitem{GS} I. Gel'fand, G. Shilov,
{\bf Generalized functions, Vol.1,} GIFML, Moscow, 1959.  (in Russian)

\bibitem{R691} A.G.Ramm, {\bf Symmetry Problems.  The Navier-Stokes Problem}, Morgan \& Claypool Publishers, San Rafael, CA, 2019.

\bibitem{R677} A.G.Ramm, Existence of the solutions to convolution equations with distributional kernels, Global Journ. of Math. Analysis, 6, N1, (2018), 1-2.

\bibitem{R684} A.G.Ramm,
Solution of the Navier-Stokes problem, Appl. Math. Lett., 87, (2019), 160-164.


\end{thebibliography}
\end{document}